\documentclass[a4,12pt]{article}

\pdfoutput=1
\usepackage{amsfonts}
\usepackage{rotating}
\usepackage{chapterbib}
\usepackage{graphicx, color}
\usepackage[latin1]{inputenc}
\usepackage[english]{babel}
\usepackage{amssymb}
\usepackage{amsmath, eurosym, layout, stmaryrd, hyperref}
\usepackage{fancyhdr}
\usepackage{times}
\usepackage{amsthm}
\usepackage{wasysym}

\theoremstyle{definition}

\renewcommand{\baselinestretch}{1.25}

\title{Invariably Suboptimal \\ \small{An attempt to improve the voting rules of Treaties of Nice and Lisbon}}

\author{Werner Kirsch and Jessica Langner}

\begin{document}
\thispagestyle{empty}
\maketitle

\renewcommand{\baselinestretch}{1.0}
\begin{abstract}\noindent
We investigate the voting rules in the Council of the European Union. It is known that both the current system, according to the Treaty of Nice, and the voting system proposed in the Lisbon treaty deviate strongly from the square root law by Penrose which is known to be the ideal voting rule under certain assumptions. In 2004 S{\l}omczynski and \.{Z}yczkowski designed a voting system, now known as the Jagiellonian Compromise, which satisfies the square root law with very high accuracy. In this system each member state obtains a voting weight proportional to the square root of the population. Then the quota is fixed in such a way that the voting power of each country is also proportional to the square root of the population.

In this paper we investigate to which extent a change of the quota in the  Treaty of Nice and the Treaty of Lisbon may bring the voting power closer to the ideal square root distribution.
Our computations show that even with optimal quota both systems are way off the ideal power distribution.
\end{abstract}
\renewcommand{\baselinestretch}{1.25}


\section{Introduction}\label{sec:1}


Political power, notably voting power, can be measured through the concept of Power Indices \cite{FM98}, \cite{TP08}. In particular, the Banzhaf Index measures how frequently a given voter
is decisive in a voting panel if all voting outcomes are counted with the same weight (for details see \cite{TP08} or our discussion below).

The Council of the European Union is a typical example  where Power Indices can help to understand the power structure within this legislative body. In fact, with its current 27 members
and the complicated decision rules voting in the Council is hard to analyze without using mathematical tools. The Council consists of one representative of each member state. The
members of the Council have different voting power  depending, in a nonsystematic way, on the size of the country they represent. The current voting system, according to the Treaty of Nice, has three
components. The first component requires the majority of states, the second a qualified majority with respect to voting weights (see table \ref{tab:0}) assigned to the states by the treaty.
The third component requires that the supporters of a proposal represent at least 62\% of the Union's population.

The draft constitution contained a new voting system for the Council, the \lq double majority\rq. This system was adopted by the Reform Treaty (Treaty of Lisbon), which is currently under discussion,
despite its rejection in a referendum in Ireland. The double majority, as
the name suggests, has two components. To make a proposal pass the Council members supporting it must represent both a qualified majority (55\%) of the states and a qualified
majority (65\%) of the population of the European Union.

Obviously, these two voting systems are very different. In fact, they lead to very different distributions of power among the states. Thus, the question arises, how a
\emph{fair} voting system should look like. An answer to this was given by Lionel Penrose as early as 1946 \cite{Pe46}. Penrose computed, what is now known as the Banzhaf Index $\beta(N)$
for a voter in a country with population $N$. He found that $\beta(N)$ is proportional to $\frac{1}{\sqrt{N}}$. Consequently, in a multinational body, such as the Council of
the EU, with one representative per state each state should have a power proportional to the square root of the country's population. This result is known as the square root law by
Penrose \cite{Pe46}, \cite{FM98}. There are various considerations about the Council of the EU in connection with the square root law \cite{Ad05}, \cite{BBGW00}, \cite{Bi04}, \cite{Bo04}, \cite{BW04}, \cite{FM00}, \cite{HM02}, \cite{Ho00}, \cite{HT06}, \cite{Ki01}, \cite{Ko05},
 \cite{Le02}, \cite{Li04}, \cite{LM03}, \cite{LM04}, \cite{LW98}, \cite{MN07}, \cite{Mo02}, \cite{Pa05}, \cite{Pl04}, \cite{So04}, \cite{Su00}, \cite{SZ04}, \cite{SZ07a}, \cite{TW00},
 \cite{Wi03}. In particular, it is well known that both the voting rules of the
Nice Treaty and those of the Treaty of Lisbon deviate strongly from the square root distribution of power. Consequently, those voting systems distribute the voting power unequally
among the citizens of the member states \cite{ABF04}, \cite{BJ04}, \cite{BW03a}, \cite{BW03b}, \cite{BW04},
 \cite{FM04a}, \cite{FM04b}, \cite{FPS03}, \cite{Le02}, \cite{Pl03}, \cite{Pl04}, \cite{PS03}.

In 2004 two Polish scientists, Wojciech S{\l}omczynski and Karol \.{Z}yczkowski devised a voting system known as the \lq Jagiellonian Compromise\rq \cite{SZ04}, \cite{SZ07a}.
In this system each member state obtains a \emph{voting weight} proportional to the square root of its population. This does not automatically give a distribution of \emph{power}
according to the square root law. However, S{\l}omczynski and \.{Z}yczkowski observed that this is the case with a particular choice of the \emph{quota}, i.e. the threshold to
reach a qualified majority. In fact, they found that with a quota of 61.4\% the \emph{voting power} (as measured by the Banzhaf Index) agrees to a very high degree of accuracy
with the square root law.

The Jagiellonian Compromise was put forward by the Polish government on the EU summit in Brussels in 2007. However, the heads of states and governments rejected this system in favor of
the double majority. Presumably, at this late stage the summit did not want to change the voting rules completely. One might hope that it would be much easier to keep either the
basic rules of the Nice Treaty or those of the Lisbon Treaty and modify a few voting rules. For example, one could just change the quota involved in such a way that one gets closer to the
square root law.

It is our task in this paper to explore to which extent one can approximate the square root distribution of power by adjusting the quota in the Nice system and for the double majority.
We compute the Banzhaf Indices for a large variety of quota for the different components of the voting systems. These results are compared to the square root law. As a measure of deviation
 from the square root law we consider the sum over all member states of the squared deviations as well as the maximum (over the states) of the deviation from the square root law.

Besides the distribution of power within the Council we also take into consideration the ability of the body to make decisions, i.e. the efficiency of the system. This value, also known as decision probability, is given by the percentage of the constellation of votes, which make a proposal pass: The higher the efficiency the easier to change the status quo, the lower the efficiency the easier to block a change. It is clear that an increase of quota will decrease the systems efficiency. While one might argue that the efficiency of the Council should not be too high to avoid
domination of a big minority of states by a small majority, the efficiency must also be not too low to ensure the EU's ability to make decisions at all.

This paper is organized as follows. In section \ref{sec:2} we give a brief introduction of the voting systems towards the Treaty of Nice and the Treaty of Lisbon. In section \ref{sec:2} we introduce the theory of voting power and a fair distribution of voting weights. In this context the square root law of Penrose will be explained. The Jagiellonian Compromise, a voting system which fulfills the square root law, is mentioned in section \ref{sec:3}. With our acquired knowledge we analyze the two treaties in detail, in particular, the obvious defects concerning the distribution of voting weights, voting power and the effectiveness will be discussed in section \ref{sec:4}. The fifth section is the main part of this paper. Here, we introduce our course of action to improve the two treaties towards the principle of equality under European citizens. We present and discuss our results and give a compromise solution for the current state of affairs. The last section of this paper contains concluding remarks.

\section{Two voting systems for the Council}
Since 2001 decision rules for voting in the European Union Council of Ministers are laid down in the Treaty of Nice: Each Member State of the European Union is assigned a voting weight (see table 1) which is a result of negotiation among the Member States. This value reflects to some degree the country's population. The Council adopts a proposal if the following three conditions (``triple majority'') are satisfied:
\begin{enumerate}
  \item The sum of the weights of the Member States vote in favor is at least 255 (of 345).
  \item A simple majority of Member States vote in favor is required (14 of 27).
  \item The Member States forming the simple majority represent at least 62\% of the overall population of the European Union.
\end{enumerate}

\begin{table}[h]
\begin{center}
\caption[]{} \label{tab:0}
\begin{tabular}{|l|c||l|c||l|c|} \hline
\multicolumn{6}{|c|}{Voting weights by the Treaty of Nice} \\\hline
Member State & Weight & Member State & Weight & Member State & Weight \\\hline\hline
Germany & 29 & Belgium & 12 & Finland & 7  \\\hline
France & 29 & Portugal & 12 & Ireland &  7  \\\hline
United Kingdom & 29 & Czech Republic &  12 & Lituania &  7  \\\hline
Italy & 29 & Hungary &  12 & Latvia &  4  \\\hline
Spain &  27 & Sweden &  10 & Slovenia &  4  \\\hline
Poland & 27 & Austria &  10 &  Estonia &  4  \\\hline
Romania &  14 & Bulgaria &  10 & Cyprus &  4  \\\hline
Netherlands &  13 & Denmark &  7 &  Luxembourg &  4  \\\hline
Greece &  12 & Slovak Republic & 7 & Malta &  3 \\\hline
\end{tabular}\end{center}
\end{table}

Mathematical analyses have shown that the three voting criteria have different effects as far as the voting outcome concerned. The first condition is the most significant one: If a qualified majority of voting weights is achieved, then in the most instances these voting weights are given by a simple majority of Member States. In contrast, the third condition has a much similar effect on the voting outcome: The probability of forming a coalition which would meet only the first and second but not the third condition is extremely low \cite{FM01}, \cite{Ki01}. Moreover, most experts agree that Nice has major drawbacks. A first one lies in the decision making efficiency of the voting body. The decision making efficiency is equal to the probability that a random proposal will be passed by vote. Here, the value of this quantity is very low with 2.03\%\footnote{Data for calculations are used from EUROSTAT: First results of the demographic data collection for 2008 in Europe.}. There already exist publications about modifying Nice such that voting power doesn't change fundamentally but its formal effectiveness increases significantly \cite{BBGW00}, \cite{BW04}. A second drawback of Nice lies in the required efforts to extend the European Union. Any extension of the Union needs a new negotiations of voting weights and thresholds.

From 2014 on an alternative voting system laid down in the Treaty of Lisbon should replace the current voting system according to the Treaty of Nice. The Treaty of Lisbon was signed in Rome in 2004 but it is not ratified until now, i.e. the Republic of Ireland has rejected it. According to the Treaty of Lisbon the Council adopts a proposal if the following two criteria (``double majority'') are satisfied:
\begin{enumerate}
  \item At least 55\% of the Member States vote in favor is required (15 of 27).
  \item The Member States forming the qualified majority represent at least 65\% of the overall population of the European Union.
\end{enumerate}
In addition, a blocking minority must include at least four Members, failing which the qualified majority shall be deemed attained. We disregard this last condition because it has no appreciable effect. The same procedure is also contained in the draft constitution of the European Convention.

The voting system for the Council according to the Treaty of Lisbon is less complex than the current system of Nice because only two criteria must be satisfied. More precisely, there are no extra weights for each state like appointed in condition 1 of the Nice Treaty. Voting weights according to Lisbon are applied directly proportional to the population of each individual Member State and the decision making efficiency is reasonably balanced with a value of 12.83\%. Moreover, any further extension of the Union is easy practicable because there is an explicit procedure how to calculate the voting weights.

Summarizing, one might receive the impression that the voting system according to Lisbon is ``better'' or ``more fair'' than the one according to Nice. Analyses have shown that this is not the case: A fair voting system of the European Union Council of Ministers should be based on a compromise between the two principles: ``equality of Member States'' and, in particular, ``equality of citizens''. Both the Treaty of Nice and the Treaty of Lisbon violate these two fundamental requirements. We will verify this statement due to concepts of the theory of voting power and its fair distribution.

\section{The theory of voting power}\label{sec:2}
Voting systems consist of a set of voters and voting rules. The voting rules determine whether a proposal is accepted or not. Frequently, there are voting weights assigned to each voter. Additionally a decision threshold is defined: a proposal will be passed if the sum of the weights of the members, who vote in favor, meets or exceeds the given threshold.

An important aspect of voting systems is the political power of the members which is also known as \emph{voting power}. Voting power is a mathematical concept which quantifies the influence a voter has on election at the system. Its theory can be traced back to works of Penrose and Banzhaf \cite{Pe46}, \cite{Ba65}. (See also \cite{SS54}, \cite{DP78}, \cite{Jo78} for alternative concepts.) Assume a member can either vote in favor or against a proposal within a decision. Then he or she has influence on the decision if he or she can turn the voting outcome by changing his or her voting behavior (to make the proposal pass by voting in favor and to make it fail otherwise). In such a situation a member is decisive. This \emph{decisiveness} is the basic idea behind voting power \cite{FM98}.

There are several methods to measure the voting power of a member. These methods are developed in the theory of the indices of political power (see books \cite{FM98}, \cite{TP08}). Power indices count in different ways in what extend an actor is decisive. One of the most popular ones is the Banzhaf Index \cite{Ba65}. The Banzhaf Index measures the a priori voting power of each member of a voting body without any previous knowledge of the single voters. Therefore it is natural to assume that all potential coalitions are equally likely. In the course of an extension of the European Union a priori power indices are useful to apply the mentioned principles of equality equally to new members as well as old ones (see also \cite{Le03}). Thus, we use the Banzhaf Index for analyses of distribution of voting power of the members in the Council \cite{ABF04}, \cite{BBGW00}, \cite{BJ04}, \cite{BW03a}, \cite{BW03b}, \cite{BW04}, \cite{FM00}, \cite{FM04a}, \cite{Le02}, \cite{Pl04}, \cite{SZ04}.\cite{SZ07a}. The Banzhaf Index is defined as follows. Assume $n$ is the number of members of a voting system. Consider each possible coalition within a member $i$. These are $2^{n-1}$. Then, the total Banzhaf Index of $i$, $TB_i$, is equal to the number of coalitions for which $i$ is decisive. The normalized Banzhaf Index of $i$, $NB_i$, is equal to the probability that $i$ is decisive: $NB_i=\frac{TB_i}{2^{n-1}}$. Finally, the percentage of influence $i$ has is given by the Banzhaf Index of $i$, $\beta_i=\frac{TB_i}{\sum_{j=1}^nNB_j}$. This quantity expresses the relative share of potential voting power of a member $i$ in the voting body. For example see the distribution of voting power of the European Economic Community of 1958-1972 in \cite{TP08}. Additionally, little shifts of quota yield to different voting power distributions.

Generally, the voting power of a member is not equal to his voting weight. This is due to the situation that voting power held by a given country depends not only on its voting weight but also on the distribution of the weights among all remaining Member States. In the case of the Council voting power should be distributed primarily satisfying equality under European citizens. A citizen has influence on an election in his country only if the other voters are split in two equal parts if a simple majority of votes in favor is required. The probability that this happens is approximately proportional to the inverse of the square root of the number of citizens (see for examples \cite{Ki04}, \cite{KMSZ04}, \cite{SZ04}, \cite{SZ06}, \cite{SZ07a}). So, if a country has $N$ citizens, then the influence of a citizen on a country's decision is proportional to $\frac{1}{\sqrt{N}}$. If we want to give all citizens the same influence on the Council's decision regardless of their home country we have to assign voting power in the Council proportional to $\sqrt{N}$. This is the square root law of Penrose.

Summarizing, the ideal distribution of voting power in the case of indirect voting consists of the Banzhaf Indices $\beta_{0i}=\frac{\sqrt{N_i}}{\sum_{j=1}^n\sqrt{N_j}}$ for each member $i$. Here $N_i$ represents the population factor of each state $i$.

To obtain a system with an ideal or fair distribution of voting power it is obvious to choose the voting weight of each Member State proportional to the square root of its population, thus equal to $\sqrt{N_i}$. This is neither a necessary nor a sufficient condition. Finally, the distribution of voting power depends on the quota (the threshold to make a proposal pass). There exists an optimal quota for which the voting power of any state is proportional to its voting weight \cite{FLMR07}, \cite{SZ06}, \cite{SZZ06}, \cite{SZ07b}. To gain this optimal quota $q_0$ we use the method of least squares: That choice of $q$ which has its least value of the sum of squared residuals $\sigma_q$ is our demanded quota. Thus, we minimize the value of the term $\sigma_q^2=\sum_{i=1}^n(\beta_{0i}-\beta_{qi})^2$ which depends on the given quota $q$. $\sigma_q^2$ is also called \emph{error rate}. In addition, the value of $\frac{\beta_{0i}-\beta_{qi}}{\beta_{0i}}$ expresses the relative deviation between demanded and obtained voting power. In the case of a minimal error rate voting weights and voting power equals best possibly. The less the error rate $\sigma_{q_0}$ the more transparent the system.

In 2007 S{\l}omczynski and \.{Z}yczkowski presented a simple mathematical formula to approximate such a quota $q_0$, in particular $q_0=\frac{1}{2}\left(1+\frac{\sqrt{N_1+\dots+N_n}}{\sqrt{N_1}+\dots+\sqrt{N_n}}\right)$, which yields to $q_0=61.57\%$\footnote{Data from EUROSTAT: First results of the demographic data collection for 2008 in Europe.} \cite{SZ07b}. This new calculation method can only be used if there are voting weights distributed as proposed by Penrose and equality under citizens with only one voting criterion is required. If there is a different voting weight distribution or there exist more than one quota to determine we have to work with the least squares method.

The European Union is not only a union of individuals but also a union of states. An additional requirement of a simple majority of Member States (``One State, One Vote'') would cause only a moderate deviation from the ideal case \cite{Ki04}, \cite{KMSZ04}, \cite{KSZ07}. Indeed a new optimal quota $q^\ast_0$ can be calculated with less discrepancies in the voting power distribution than with the previous $q_0$. Therefore we have to use the method of least squares again.

Beyond a fair distribution of influence we should consider the effectiveness of a system. Effectiveness is equal to the decision probability the voting body passes a proposal. This quantity is also called the Coleman power of a collectivity act \cite{Co71}. Assuming that all coalitions are equally likely its value is given by the percentage of the constellation of votes, which make a proposal pass: The higher the effectiveness the easier to change the status quo, the lower the rate the easier to block a change. So, the degree of the effectiveness depends on the given voting rules, in particular the quotas.

Voting systems based on the square root low of Penrose were proposed and discussed many times. One of the best-known proposals is the Jagiellonian Compromise.

\section{The Jagiellonian Compromise}\label{sec:3}
In 2004 the polish scientists, Wojciech S{\l}omczynski and Karol \.{Z}yczkowski, from the Jagiellonian University of Kraków, Poland, presented a voting system for the Council of Ministers of the European Union, the Jagiellonian Compromise \cite{SZ04}, \cite{SZ06}, \cite{SZ07a}. They constructed a voting system as follows: The voting weight of each Member State is chosen according to the square root law of Penrose, thus equal to $\sqrt{N_i}$ where $N_i$ is the population factor of the $i$-th Member State. Then, an optimal quota $q$ is calculated using the methods above. The Jagiellonian Compromise is also known as $P-q\%$ solution due to the work of Penrose.

With current population data we gain an optimal quota $q_0=61.5\%$ with a minimal error rate of $0.00005$\permil.  Our analyses have shown that the maximal relative deviation between $\beta_{0i}$ and its corresponding $\beta_{q_0i}$ is about the less value of $0.14\%$. In addition, the effectiveness value is about $16.43\%$. For voting weights and voting power see table 2.

\begin{table}[h]
\begin{center}
\label{Tab:2}\caption[]{}
\begin{tabular}{|l|r|r|c|c|} \hline
\multicolumn{5}{|c|}{Distribution of votes and voting power on the Council of Ministers} \\\hline\hline
\multicolumn{5}{|c|}{The Jagiellonian Compromise - P-61.5 solution}\\\hline
Member State & Population & Population & voting weight& Banzhaf Index \\
 &  & square root & in \% & $\beta_i$ in \% \\\hline\hline
Germany & 82.221.808 & 9.067,6242 & 9,4108 & 9,3978 \\\hline
France & 63.753.140 & 7.984,5563 & 8,2867 & 8,2933 \\\hline
United Kingdom & 61.185.981 & 7.822,1468 & 8,1181 & 8,1254 \\\hline
Italy & 59.618.114 & 7.721,2767 & 8,0135 & 8,0214 \\\hline
Spain & 45.283.259 & 6.729,2837 & 6,9839 & 6,9924 \\\hline
Poland & 38.115.641 & 6.173,7866 & 6,4074 & 6,4141 \\\hline
Romania & 21.528.627 & 4.639,8951 & 4,8155 & 4,8175 \\\hline
Netherlands & 16.404.282 & 4.050,2200 & 4,2035 & 4,2038 \\\hline
Greece & 11.214.992 & 3.348,8792 & 3,4756 & 3,4746 \\\hline
Belgium & 10.666.866 & 3.266,0168 & 3,3896 & 3,3885 \\\hline
Portugal & 10.617.575 & 3.258,4621 & 3,3818 & 3,3807 \\\hline
Czech Republic & 10.381.130 & 3.221,9761 & 3,3439 & 3,3428 \\\hline
Hungary & 10.045.000 & 3.169,3848 & 3,2893 & 3,2881 \\\hline
Sweden & 9.182.927 & 3.030,3345 & 3,1450 & 3,1437 \\\hline
Austria & 8.331.930 & 2.886,5083 & 2,9957 & 2,9939 \\\hline
Bulgaria & 7.640.238 & 2.764,098 & 2,8687 & 2,8671 \\\hline
Denmark & 5.475.791 & 2.340,0408 & 2,4286 & 2,4269 \\\hline
Slovak Republic & 5.400.998 & 2.324,0047 & 2,4119 & 2,4100 \\\hline
Finland & 5.300.484 & 2.302,278 & 2,3894 & 2,3876 \\\hline
Ireland & 4.419.859 & 2.102,3461 & 2,1819 & 2,1801 \\\hline
Lituania & 3.366.357 & 1.834,7635 & 1,9042 & 1,9025 \\\hline
Latvia & 2.270.894 & 1.506,9486 & 1,5640 & 1,5623 \\\hline
Slovenia & 2.025.866 & 1.423,3292 & 1,4772 & 1,4757 \\\hline
Estonia & 1.340.935 & 1.157,9875 & 1,2018 & 1,2003 \\\hline
Cyprus & 794.580 & 891,3922 & 0,9251 & 0,9241 \\\hline
Luxembourg & 483.799 & 695,5566 & 0,7219 & 0,7210 \\\hline
Malta & 410.584 & 640,7683 & 0,6650 & 0,6642 \\\hline\hline
Sum & 497.481.657 & 96353,8647 & 100 & 100 \\
  \hline
\end{tabular}\end{center}
\end{table}

Some advantages arise from the proposed voting system: First of all it is \emph{simple}, because it is based on a single criterion, more precisely, only one condition must be satisfied. It is \emph{neutral} by reason that it cannot a priori favor or handicap any Member State. It is \emph{fair}, because every citizen has the same potential influence on decisions regardless from his home country. It is \emph{transparent} in the case that voting power and voting weight are almost equal. It is \emph{easy extendible}: any new Member State achieves a voting weight proportional to the square root of its population factor. Solely a new optimal quota must be calculated. It is \emph{moderately efficient} because an addition of Member States does not decrease the effectiveness.

On closer observation the additional requirement of a simple majority of Member States (in the following denoted by JC+) is postulated with an equal $q_0$. That yields to an error rate of $0.07425\permil$. The relative voting power deviation takes a maximum value of $30.64\%$. This is no more a moderate deviation from the ideal case. Only the effectiveness value almost levels off with $16.08\%$. Observing the least squares we gain a new optimal quota of $q^\ast_0=64.7\%$. Here, the error rate takes its minimum value of $0.03275\permil$. This is only the half of the error rate value \clearpage than with an unchanged quota. The maximal relative deviation is only about $11.68\%$. This is nearly one third in comparison with $q_0$. However, the effectiveness decreases on the lower value of $10.39\%$. In terms of an as best as possible fair distribution of voting the quota $q^\ast$ should be applied. Figure 1 shows that the voting power of the JC+ 64.7-solution is better approximated to Penrose's $\beta_0$ than the JC+ 61.5-solution.


\begin{figure}[h]
\begin{center}
\label{fig:1}\caption[]{}
\includegraphics[width=14cm]{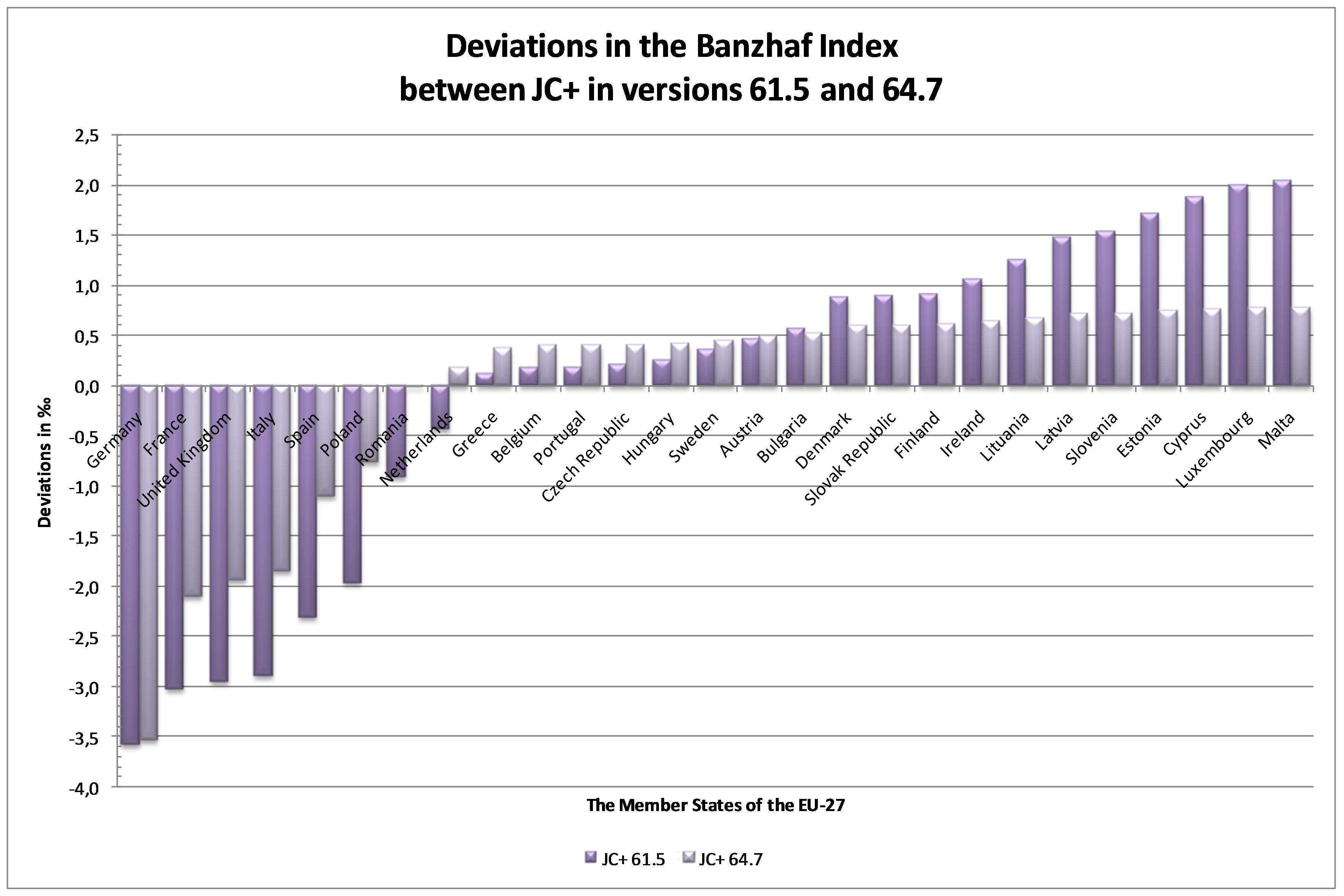}
\end{center}
\end{figure}


Various efforts have been made to promote this ideal system. For examples, in 2004 about 50 scientists supported the Jagiellonian Compromise in an open letter to the governments of the European Union Member States with the title ``Scientists for a democratic Europe''. Moreover, in the course of the EU summit in Brussels in 2007 the polish mission statement ``The square root or death'' made the problem the subject of discussions again, unfortunately, without observable success.

\section{Penrose vs. the Treaty of Nice and the Treaty of Lisbon}\label{sec:4}
With the acquired knowledge about voting power and its fair distribution we will have a second look on the two treaties for the European Union. The voting power values according to the Treaty of Nice and the Treaty of Lisbon are shown in table 3. First of all, both systems are \emph{not simple}, because at least two conditions must be satisfied. The votes, weights and thresholds for the Council laid down in both treaties are not scientific based chosen. They are results of negotiations among the Member States. Thus, both systems are \emph{not objective}. In addition, they violate the square root law, because the voting power is not distributed like Penrose specified.

This is obviously shown by the corresponding error rates and deviations in voting power (see figures 2-5). According to the Treaty of Nice the four biggest states, Germany, France, United Kingdom and Italy, are assigned with too little voting power in comparison with the optimal voting weights $\beta_0$. In contrast, especially Spain and Poland achieve too much power also some middle Member States and the four smallest Estonia, Cyprus, Luxembourg and Malta. Here, the maximal relative deviation in voting power is about $73.18\%$. The corresponding error rate has a value of $0.6052\permil$.

According to the Treaty of Lisbon especially the bigger and the smaller Member States achieve more voting power. Thus the middle States are suffered of this. Here the error rate is very high with $1.2438\permil$. Also the maximal relative deviation is about $137.53\%$. Hence, both systems are \emph{not representative}. They do not fulfill the principle ``One Person, One Vote''.  Moreover, it is quite obvious that voting weight and resulting voting power deviate strongly in comparison to the deviations in obtained and demanded voting power according to the Jagiellonian Compromise. Thus, both systems are \emph{not transparent}. According to the Treaty of Lisbon an extension of the European Union is easy practicable in the way that it needs no negotiations among the Member States about voting weights. Contrariwise, the Treaty of Nice always required new debates.

As denoted above the decision making efficiency according to the Treaty of Nice is very low: It is about $2.03\%$. The effectiveness of the Treaty of Lisbon is about $12.83\%$ which is only a little bit higher than in the case of the Jagiellonian Compromise with an additional requirement of a simple majority of Member States.

Summarizing, Nice has less power distortions than Lisbon but it is more complex. Unfortunately, the effectiveness of Nice hardly allows changes of the status quo. The current voting rules, weights and thresholds were fully discussed whole nights long without scientific based background. In addition, we know from the European Economic Community of 1958-1973 that little shifts of the quota yield to different voting power distributions. It would be a fortunate coincidence if the current thresholds produce the best possible error rate measured by the ideal voting weights due to Penrose.

\begin{table}[h]
\begin{center}
\label{tab:3}\caption[]{}
\begin{tabular}{|l||r|c|c||c|c|c|} \hline
\multicolumn{7}{|c|}{Distribution of votes and voting power on the Council of Ministers} \\\hline\hline
& \multicolumn{3}{|c||}{The Treaty of Lisbon 55/65}& \multicolumn{3}{|c|}{The Treaty of Nice 255}\\\hline
Member State & Population & Population & Banzhaf Index & Votes in & Weight & Banzhaf Index\\
& & in \% & $\beta_{i}$ in \% & the Council & in \% & $\beta_{i}$ in \% \\\hline\hline
Germany & 82.221.808 & 16,53 & 11,5362 & 29 & 8,41 & 7,7828 \\\hline
France & 63.753.140 & 12,82 & 9,0667 & 29 & 8,41 & 7,7828 \\\hline
United Kingdom & 61.185.981 & 12,30 & 8,7322 & 29 & 8,41 & 7,7827 \\\hline
Italy & 59.618.114 & 11,98 & 8,5360 & 29 & 8,41 & 7,7827 \\\hline
Spain & 45.283.259 & 9,10 & 6,6893 & 27 & 7,83 & 7,4199 \\\hline
Poland & 38.115.641 & 7,66 & 5,6050 & 27 & 7,83 & 7,4198 \\\hline
Romania & 21.528.627 & 4,33 & 4,1306 & 14 & 4,06 & 4,2591 \\\hline
Netherlands & 16.404.282 & 3,30 & 3,4952 & 13 & 3,77 & 3,974 \\\hline
Greece & 11.214.992 & 2,25 & 2,8747 & 12 & 3,48 & 3,6843 \\\hline
Belgium & 10.666.866 & 2,14 & 2,8092 & 12 & 3,48 & 3,6843 \\\hline
Portugal & 10.617.575 & 2,13 & 2,8033 & 12 & 3,48 & 3,6843 \\\hline
Czech Republic & 10.381.130 & 2,09 & 2,7750 & 12 & 3,48 & 3,6843 \\\hline
Hungary & 10.045.000 & 2,02 & 2,7349 & 12 & 3,48 & 3,6843 \\\hline
Sweden & 9.182.927 & 1,85 & 2,6321 & 10 & 2,90 & 3,0924 \\\hline
Austria & 8.331.930 & 1,67 & 2,5302 & 10 & 2,90 & 3,0924 \\\hline
Bulgaria & 7.640.238 & 1,54 & 2,4478 & 10 & 2,90 & 3,0924 \\\hline
Denmark & 5.475.791 & 1,10 & 2,1891 & 7 & 2,03 & 2,1809 \\\hline
Slovak Republic & 5.400.998 & 1,09 & 2,1803 & 7 & 2,03 & 2,1809 \\\hline
Finland & 5.300.484 & 1,07 & 2,1681 & 7 & 2,03 & 2,1809 \\\hline
Ireland & 4.419.859 & 0,89 & 2,0625 & 7 & 2,03 & 2,1809 \\\hline
Lituania & 3.366.357 & 0,68 & 1,9362 & 7 & 2,03 & 2,1809 \\\hline
Latvia & 2.270.894 & 0,46 & 1,8044 & 4 & 1,16 & 1,2502 \\\hline
Slovenia & 2.025.866 & 0,41 & 1,7747 & 4 & 1,16 & 1,2502 \\\hline
Estonia & 1.340.935 & 0,27 & 1,6920 & 4 & 1,16 & 1,2502 \\\hline
Cyprus & 794.580 & 0,16 & 1,6260 & 4 & 1,16 & 1,2502 \\\hline
Luxembourg & 483.799 & 0,10 & 1,5886 & 4 & 1,16 & 1,2502 \\\hline
Malta & 410.584 & 0,08 & 1,5796 & 3 & 0,87 & 0,9422 \\\hline\hline
Sum & 497.481.657 & 100,00 & 100,0000 & 345 & 100 & 100 \\\hline
\end{tabular}\end{center}
\end{table}
\clearpage


\section{Improvements}\label{sec:5}
It is our goal to optimize the current Treaty of Nice and the Treaty of Lisbon. Therefore, we fractionally modify the voting rules: The existing voting weights will be unchanged retained and only the several (up to three) thresholds will be shifted. We search for a constellation of quotas such that the resulting Banzhaf Indices reach the least possible error rate. Therefore, we have programmed a Java-applet which calculates for several thresholds tuples the Banzhaf Index values of each Member State, the corresponding error rates, the maximal deviations between demanded and obtained voting power and the effectiveness of the voting systems. As basis data we use the population values from EUROSTAT 2008 (see table 2).

The Treaty of Nice will be investigated with an unchanged simple majority of Member States, thus 14. The quota of the sum of voting weights (currently 255 (=73,91\%)) will be shifted from 190 (=55.07\%) up to 275 (=79.71\%) in integers. For each given quota of voting weights we shift the overall population quota (currently 62\%) from 51\% up to 85\% in steps of 1\%.

The Treaty of Lisbon will be analyzed with integer majority of Member States from 14 up to 18 (currently 15). A majority of 14 states relates to a relative majority of 48.15\% up to 51.85\% $(\frac{13}{27}=48.15\%)$, 15 up to 55.55\%, also 18 up to 66.66\%. For each given integer majority, we shift the overall population quota (currently 65\%) from 51\% up to 85\% in steps of 0.1\%.

Beyond our boundary values the error rate significantly increases. This is due to the fact, that a higher quota give more power to smaller states (a proposal will be passed with almost unanimity) and lower quota more power to bigger states. Furthermore, we want to include the corresponding effectiveness value within our approach of optimization. It is easy to see that the decision making efficiency goes to zero with increasing quota.

In the case of the Treaty of Nice our calculations have produced the threshold tuple (14/263/80\%) due to the least minimal error rate of $0.2286\permil$. Compared with the Jagiellonian Compromise ($0.00005\permil$) Nice's best possible error rate still deviates strongly from the ideal case. This is also indicated by a maximal relative deviation in voting power with $42.9\%$ (JC: 0.14\%). Therewith, the effectiveness is very low with $0.99\%$.

In the case of the Treaty of Lisbon our calculations have produced the threshold tuple (17/77.5\%) due to the least minimal error rate of $0.52118\permil$. Compared with Nice the best possible error rate of Lisbon is additionally 127\% higher. This is also indicated by a maximal relative deviation in voting power with $135.51\%$ (JC: 0.14\%). Concluding, the effectiveness is very low with $2.23\%$ thus near to Nice in its current version.

Using the examples of Nice (14/255/.), (14/263/.) and Lisbon (15/.), (17/.) the development of error rates and effectiveness by shifting the population quota is pictured in the diagrams 2 and 3. Comparisons of Nice in versions (14/255/62\%) and (14/263/80\%) and Lisbon in versions (15/65\%) and (17/77.5\%) are shown in the figures 4 and 5. In particular, in both cases Germany's large deviation between demanded and obtained influence is strongly decreased. According to Nice both Spain and Poland are still assigned too much voting power as specified by Penrose. According to Lisbon middle-size States are still assigned with little voting power and small and big States with too much.

Summarizing, our optimized threshold tuples produce less deviations in voting power measured by the least possible error rate than the current versions. Nevertheless, each quota constellation produce a \emph{significant} deviation to Penrose's ideal case. For comparison, in the appendix we have listed several threshold tuples with fixed voting weight and State quota and optimized population quota such with minimal error rate. In addition, the values of the related maximal relative deviation in voting power show that the resulted systems are neither transparent nor representative. Moreover, these optimizations due to the error rates lead to a very low effectiveness. Thus, in such voting systems it would be easy to block proposals.

Due to these results we reconsider our analyses and include the effectiveness values in our solution approach. One approach might be to find a compromise between current error rate, optimal error rate and a reasonable decision probability. Therefore, we are geared to the effectiveness value of the Jagiellonian Compromise including the requirement of a simple majority of Member States thus $10.39\%$. According to the Treaty of Nice we refer to the threshold tuple (14/220/66\%). We gain an error rate of $1.07\permil$, a maximal relative deviation of $37.43\%$ and an effectiveness of $10.52\%$. According to Treaty of Lisbon we refer to the threshold tuple (15/67.5\%). This yield to an error rate of $1.5975\permil$, a deviation of $106.62\%$ and an effectiveness of $10.36\%$. Certainly, there are several solutions for constellations of effectiveness and error rate values supposable. But, by now it might be conceivable that it needs many new debates among the Member States and a lot of time to find such distribution keys.

\begin{figure}[h]
\begin{center}\caption[]{}
\includegraphics[width=14cm, height=9cm]{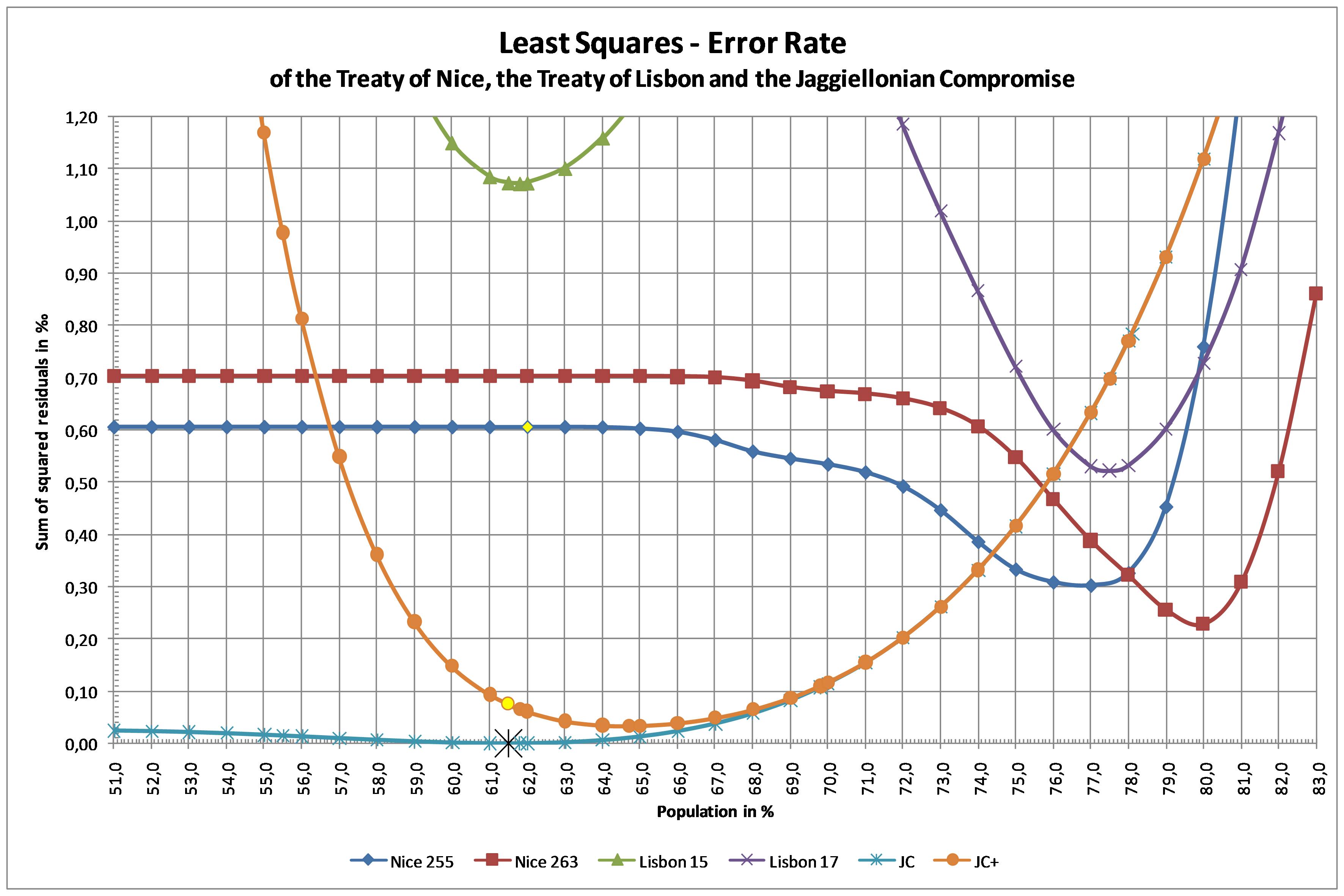}\caption[]{}
\includegraphics[width=14cm, height=9cm]{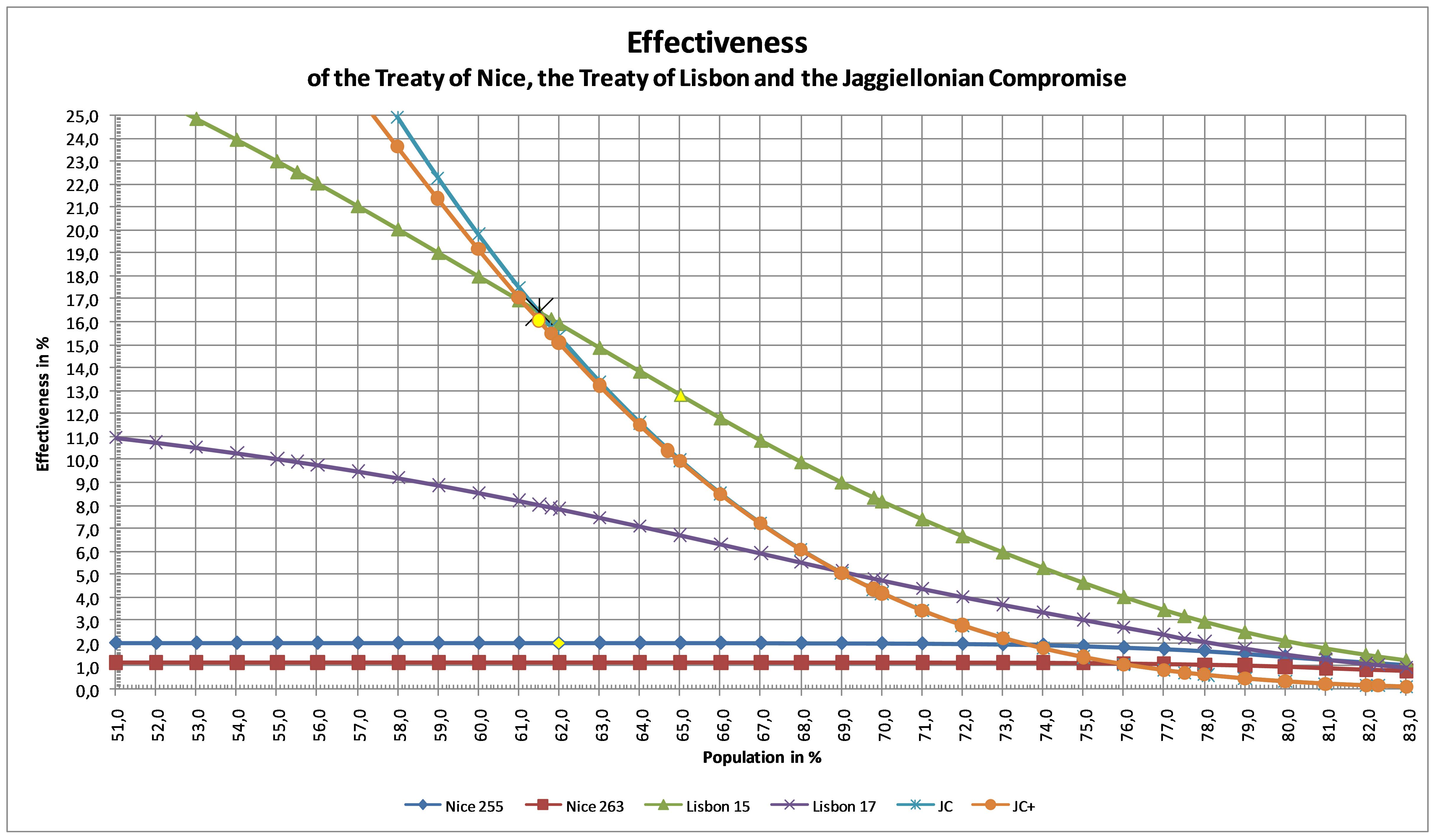}
\end{center}
\end{figure}
\clearpage
\begin{figure}[h]
\begin{center}\caption[]{}
\includegraphics[width=14cm, height=9cm]{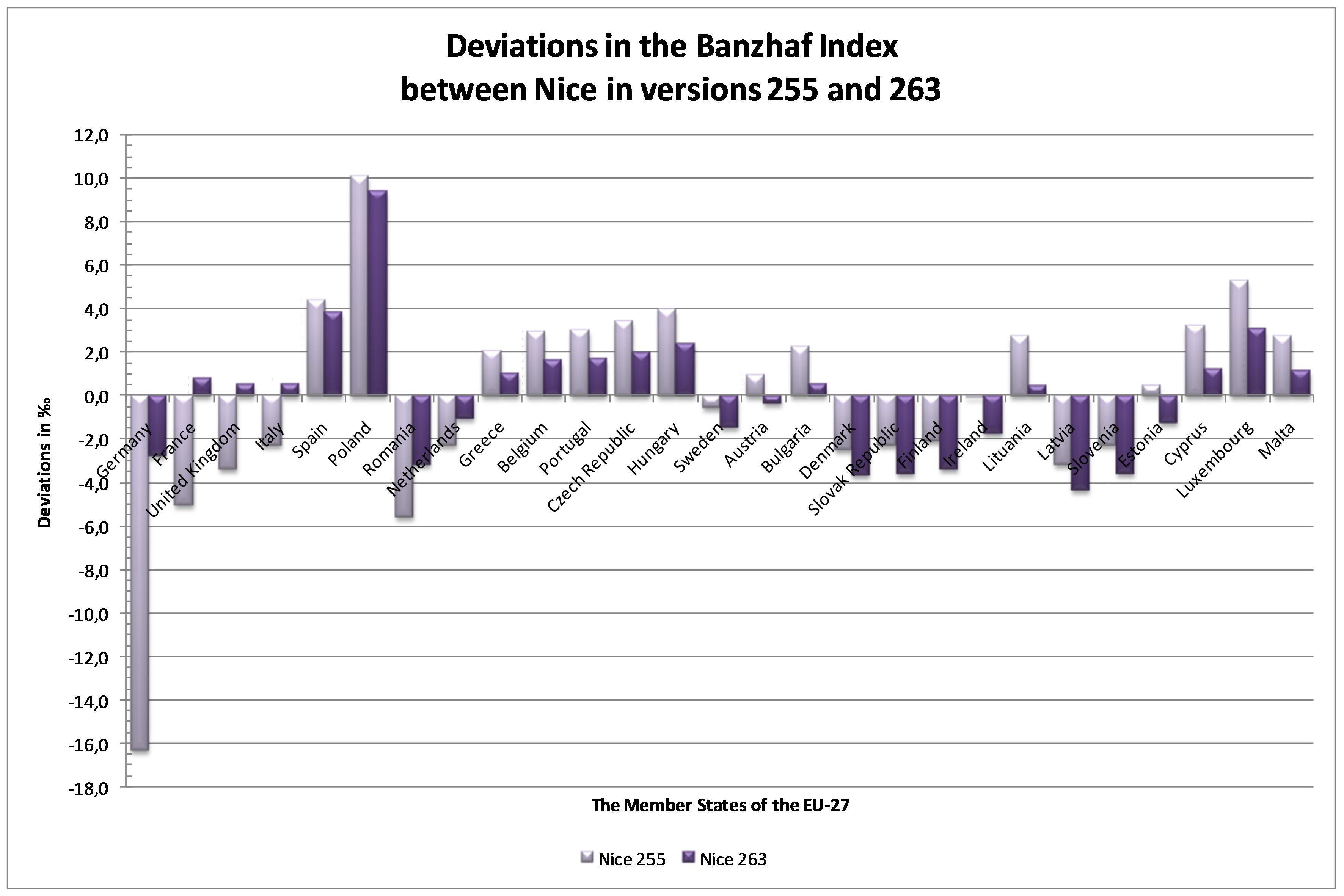}\caption[]{}
\includegraphics[width=14cm, height=9cm]{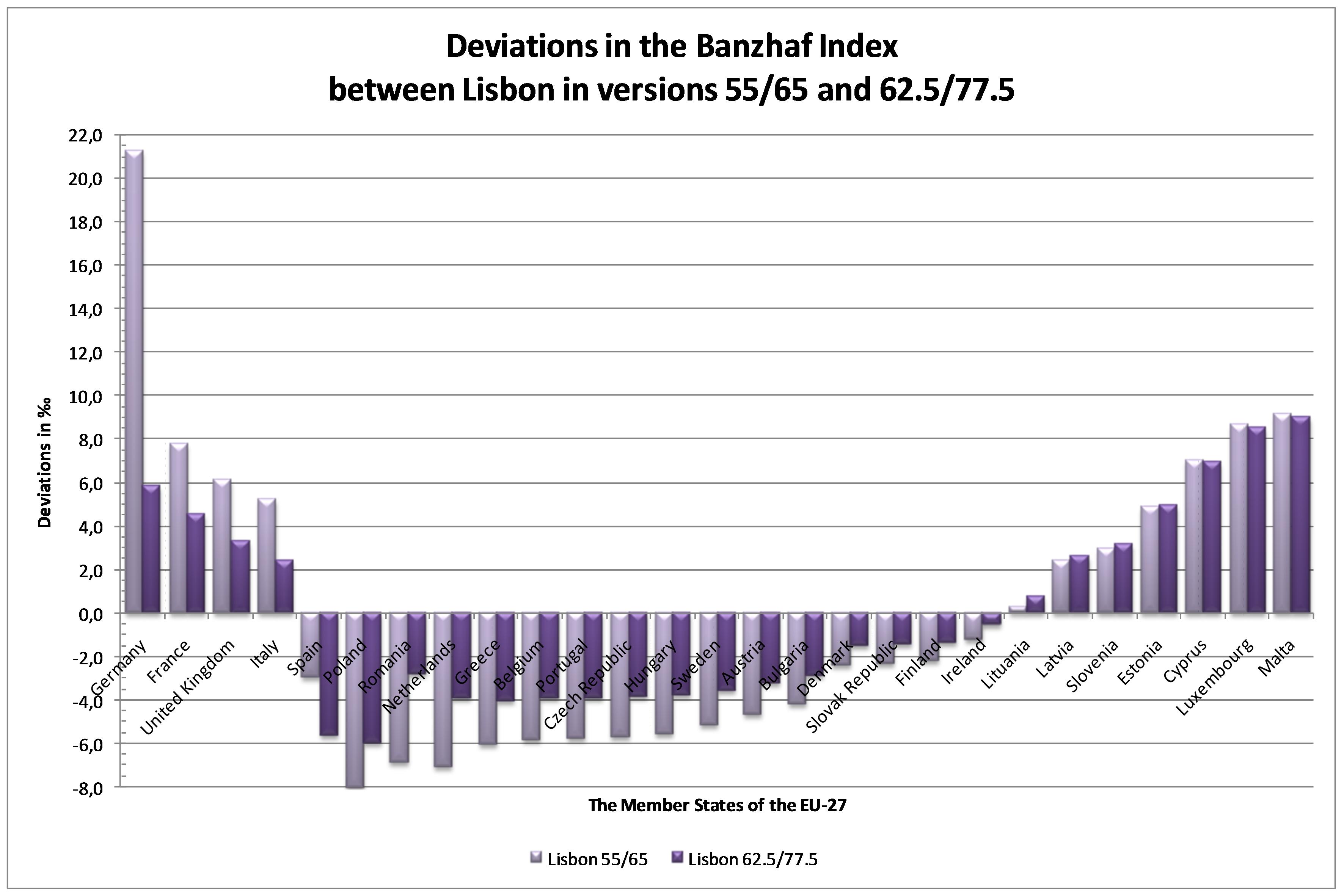}
\end{center}
\end{figure}
\clearpage

\section{Conclusions}\label{sec:5}
As several publications have shown both voting systems for the Council of Ministers of the European Union according to the Treaty of Nice and the Treaty of Lisbon deviate strongly from Penrose's solution of a fair distribution of voting power in such a voting body. In this publication we tried to improve these treaties with respect to such a fair distribution. To do so we modified the voting rules by keeping the voting weights and only shifting the thresholds. This procedure results only in a modest improvement of the system. Even with optimal quota both systems deviate strongly from a fair distribution of power. Thus, both the Treaty of Nice and the Treaty of Lisbon turn out to be invariably suboptimal.

As a consequence the voting system for the Council has to be changed in a more fundamental way than merely adjusting quota. It seems to us that the Jagiellonian Compromise is a good basis for a new voting system.



\vspace{10mm}
\emph{Werner Kirsch.}\\ Fakultät für Mathematik und Informatik, FernUniversität in Hagen, 

\vspace{2.5mm}
\emph{Jessica Langner.}\\ Fakultät für Mathematik, Ruhr-Universität Bochum, D-44780 Bochum, 

\clearpage
\section{Appendix}

\begin{table}[h]
\begin{center}\caption[]{}
\begin{tabular}{|c|c||c|c|c|} \hline
\multicolumn{5}{|c|}{The optimal threshold values for the Treaty of Lisbon} \\\hline
Quota of & Quota of & Sum of square & maximal relative & Effectiveness\\
States & Population & residuals in \permil & deviation in \% & in \% \\\hline\hline
14 & 0,555 & 1,08417 & 173,55 & 28,33 \\\hline
14 & 0,600 & 1,83319 & 116,39 & 21,89 \\\hline
15 & 0,618 & 1,07222 & 180,26 & 16,13 \\\hline
15 & 0,650 & 1,24384 & 137,53 & 12,83 \\\hline
15 & 0,675 & 1,59748 & 106,62 & 10,36 \\\hline
16 & 0,698 & 0,83352 & 161,88 & 6,74 \\\hline
17 & 0,775 & 0,52118 & 135,51 & 2,23 \\\hline
18 & 0,823 & 0,75088 & 163,24 & 0,78 \\\hline\hline
\multicolumn{2}{|c||}{optimal values}  & 0,52118 & 116,39 & 28,33 \\\hline
\end{tabular}\label{tab:5}
\end{center}
\end{table}

\begin{table}[h]
\begin{center}\caption[]{}
\begin{tabular}{|c|c||c|c|c|} \hline
\multicolumn{5}{|c|}{The optimal threshold values for the Jagiellonian Compromise} \\\hline
Quota of & Quota of & Sum of square & maximal relative & Effectiveness\\
 States &  Population & residuals in \permil & deviation in \% & in \% \\\hline\hline
 & 0,615 & 0,00005 & 0,14 & 16,43 \\\hline
14 & 0,615 & 0,07425 & 30,64 & 16,08 \\\hline
14 & 0,647 & 0,03275 & 11,68 & 10,39 \\\hline\hline
\multicolumn{2}{|c||}{optimal values} & 0,00005 & 0,14 & 16,43 \\\hline
\end{tabular}\label{tab:6}
\end{center}
\end{table}

\begin{table}[h]
\begin{center}\caption[]{}
\begin{tabular}{|c|c|c||c|c|c|} \hline
\multicolumn{6}{|c|}{The optimal threshold values for the Treaty of Nice} \\\hline
Quota of & Quota of & Quota of & Sum of square & maximal relative & Effectiveness\\
 States & Weights & Population & residuals in \permil & deviation in \% & in \% \\\hline\hline
14 & 190 & 0,54 & 0,5192 & 142,08 & 27,74 \\\hline
14 & 195 & 0,55 & 0,3966 & 118,83 & 25,27 \\\hline
14 & 200 & 0,56 & 0,3388 & 95,61 & 22,11 \\\hline
14 & 205 & 0,58 & 0,3367 & 83,65 & 19,47 \\\hline
14 & 210 & 0,59 & 0,3646 & 75,48 & 16,85 \\\hline
14 & 215 & 0,60 & 0,3913 & 68,49 & 14,33 \\\hline
14 & 220 & 0,61 & 0,4099 & 63,33 & 11,97 \\\hline
14 & 225 & 0,63 & 0,4096 & 57,66 & 9,65 \\\hline
14 & 230 & 0,64 & 0,4117 & 58,76 & 7,79 \\\hline
14 & 235 & 0,66 & 0,4143 & 58,91 & 6,11 \\\hline
14 & 240 & 0,67 & 0,4188 & 61,38 & 4,77 \\\hline
14 & 245 & 0,72 & 0,4077 & 49,53 & 3,45 \\\hline
14 & 250 & 0,74 & 0,3469 & 49,47 & 2,55 \\\hline
14 & 255 & 0,62 & 0,6052 & 73,18 & 2,03 \\\hline
14 & 255 & 0,77 & 0,3016 & 45,76 & 1,76 \\\hline
14 & 258 & 0,62 & 0,6373 & 74,58 & 1,66 \\\hline
14 & 258 & 0,78 & 0,2620 & 46,73 & 1,44 \\\hline
14 & 259 & 0,79 & 0,2515 & 40,34 & 1,30 \\\hline
14 & 260 & 0,79 & 0,2391 & 43,76 & 1,23 \\\hline
14 & 261 & 0,79 & 0,2373 & 47,85 & 1,17 \\\hline
14 & 262 & 0,80 & 0,2372 & 39,38 & 1,04 \\\hline
14 & 263 & 0,80 & 0,2286 & 42,90 & 0,99 \\\hline
14 & 264 & 0,80 & 0,2318 & 47,58 & 0,93 \\\hline
14 & 265 & 0,80 & 0,2445 & 51,09 & 0,88 \\\hline
14 & 270 & 0,82 & 0,2762 & 49,75 & 0,58 \\\hline
14 & 275 & 0,84 & 0,3587 & 61,37 & 0,38 \\\hline\hline
\multicolumn{3}{|c||}{optimal values}& 0,2286 & 39,3825 & 27,74 \\\hline
\end{tabular}\label{tab:4}
\end{center}
\end{table}

\label{LastPage}
\clearpage
\end{document}